\documentclass[11pt]{amsart}

\usepackage[T1]{fontenc}
\usepackage[utf8]{inputenc}
\usepackage{lmodern}
\usepackage{microtype}
\usepackage{amsmath,amssymb,amsthm,mathtools}
\usepackage{enumitem}
\usepackage[hidelinks]{hyperref}
\usepackage[nameinlink,capitalize,noabbrev]{cleveref}

\usepackage{mathrsfs}
\theoremstyle{plain}
\newtheorem{theorem}{Theorem}[]
\newtheorem{proposition}[theorem]{Proposition}
\newtheorem{lemma}[theorem]{Lemma}
\newtheorem{corollary}[theorem]{Corollary}

\theoremstyle{definition}

\theoremstyle{remark}
\newtheorem{remark}[theorem]{Remark}

\title[Metric moduli of Enriques surfaces]{The fundamental group of the moduli space of Ricci-flat metrics on an Enriques surface}

\author{Xikai Wei}
\address{Department of Mathematics, Karlsruhe Institute of Technology, 76131 Karlsruhe, Germany}
\email{xikai.wei@kit.edu}

\date{\today}

\begin{document}

\begin{abstract}
We prove that every connected component of the coarse moduli space of smooth Ricci-flat metrics on the Enriques manifold is simply connected. The metric Teichm\"uller space is not simply connected because its period domain contains a discriminant of real codimension two. We retain this topology in an action braid group and compute the coarse fundamental group by imposing smooth metric stabilizer relations. Bertini and Kantor involutions, together with the modular description of the Enriques period quotient and its canonical stack, reduce the calculation to a cyclic group of order four. An order-eight automorphism supplies a stabilizer whose image generates this group, and a sheet-reversing metric involution completes the proof.
\end{abstract}


\maketitle

\section{Introduction}
\label{sec:introduction}

We study the coarse moduli space of smooth Ricci-flat metrics on the Enriques manifold. Every connected component is simply connected, although the corresponding metric Teichm\"uller space is not. The difference comes from stabilizers: after combining paths in the period domain with the diffeomorphism action, the relations supplied by isometries of smooth metrics kill the entire fundamental group of the coarse quotient.

Let $E$ be the smooth four-manifold underlying a complex Enriques surface. Its canonical unramified double cover $X$ is a K3 surface. The canonical bundle of $E$ has order two, so Yau's theorem provides a unique Ricci-flat K\"ahler metric in each K\"ahler class \cite{Yau1978}. Conversely, every Ricci-flat metric on $E$ lifts to a hyperk\"ahler metric on $X$. Its period can therefore be described using the invariant and anti-invariant parts of the K3 lattice.

Let $\operatorname{Met}(E)$ denote the space of smooth Riemannian metrics on $E$ with the $C^\infty$ topology, and set
$$\mathcal{R}^{\mathrm{Ric}=0}_{\mathrm{unit}}(E):=\{g\in\operatorname{Met}(E)\mid\operatorname{Ric}_g=0,\ \operatorname{Vol}_g(E)=1\}.$$
We give this space the subspace topology and define
$$\mathcal{M}^{\mathrm{Ric}=0}_{\mathrm{unit}}(E):=\mathcal{R}^{\mathrm{Ric}=0}_{\mathrm{unit}}(E)/\operatorname{Diff}(E)$$
with the quotient topology. Fundamental groups of coarse moduli spaces always mean ordinary topological fundamental groups.

\begin{theorem}[Main theorem]
\label{thm:main}
Every connected component of $\mathcal{M}^{\mathrm{Ric}=0}_{\mathrm{unit}}(E)$ is simply connected. For each such component $\mathcal{M}_0$,
$$\pi_1(\mathcal{M}_0)=\{1\}$$
\end{theorem}

The same conclusions hold without volume normalization.

Related results include the simple connectivity of the coarse moduli space of flat metrics on a torus, proved by Tuschmann and Wiemeler \cite[Proposition~5.5(1)]{TuschmannWiemeler2022}, and Degen's theorem for the unit-volume K3 metric moduli space when orbifold singularities are included \cite{Degen2020}.

More recently, Ouyang studied degenerations of Ricci-flat K\"ahler metrics on Enriques surfaces with unit-diameter normalization and described their Gromov--Hausdorff compactification in terms of the adjoint Satake compactification \cite{Ouyang2026}. This provides a complementary perspective on the global geometry of Enriques metric moduli.

In the smooth Enriques problem, an additional difficulty is the topology of the discriminant complement.

For a compatible marking of the K3 cover, the invariant and anti-invariant lattices are
$$M\cong U(2)\oplus E_8(-2),\qquad N\cong U\oplus U(2)\oplus E_8(-2).$$
The lifted metric period splits into a positive line in $M_{\mathbb{R}}$ and a positive two-plane in $N_{\mathbb{R}}$. With the orientation conventions specified below, each metric Teichm\"uller component has the form
$$\mathcal{T}_E\cong(\mathcal{D}_M\times\mathcal{D}_N^\circ)\setminus(\Delta_{\mathrm{per}}\cup\Delta_{\mathrm{met}}).$$
Both domain factors are contractible. The period discriminant $\Delta_{\mathrm{per}}$ has real codimension two and contributes meridians to the fundamental group. The remaining metric discriminant has real codimension three, so its removal does not change the fundamental group.

Let $\Gamma_E^{\mathrm{act}}$ be the effective action induced by diffeomorphisms preserving the chosen component. The action braid group retains both the period paths and their arithmetic endpoints:
$$1\longrightarrow\pi_1(\mathcal{T}_E)\longrightarrow\mathcal{B}_E\longrightarrow\Gamma_E^{\mathrm{act}}\longrightarrow1.$$
Armstrong's theorem \cite{Armstrong1968}, applied to the universal cover, identifies the coarse fundamental group with $\mathcal{B}_E/\mathcal{N}_{\mathrm{stab}}$, where $\mathcal{N}_{\mathrm{stab}}$ is normally generated by metric stabilizer classes.

We compute this quotient by constructing explicit stabilizers of
smooth Ricci-flat metrics. Allcock's generation theorem for Bertini
involutions \cite[Theorem~1.1]{Allcock2018}, together with the
K3 metric realization theory of Farb--Looijenga
\cite[Theorem~1.2]{FarbLooijenga2024}, provides the ingredients
for showing that the corresponding stabilizer relations eliminate
the kernel of the projective period action. Kantor involutions
similarly yield smooth metric stabilizers whose relations eliminate
the divisorial inertia along the retained reflection mirrors;
compare \cite[Theorem~3.3]{Allcock2018}.

We then use Borcherds's weight-four automorphic form and its
reflective divisor \cite[Theorem~3.2]{Borcherds1996}, together
with Wang--Williams's description of the graded ring of
orthogonal modular forms \cite[Theorem~5.2]{WangWilliams2020},
to compute the residual quotient. We show that this quotient
is isomorphic to the cyclic group $\mu_4$ of fourth roots of
unity. Ohashi constructed an order-eight automorphism of an
Enriques surface whose bi-canonical representation has order
four \cite[Example~1.2]{Ohashi2015}. We prove that the
corresponding smooth metric stabilizer class maps to a
generator of $\mu_4$, thereby eliminating the residual
sheet-preserving quotient. Finally, a sheet-reversing metric
involution supplies the remaining stabilizer relation and
completes the calculation of the fundamental group.

Section~\ref{sec:metric-periods} establishes the period model and actual action. Section~\ref{sec:action-braid-groups} develops the required topology. Sections~\ref{sec:metric-stabilizers} and~\ref{sec:residual-period-group} construct the stabilizers and compute the residual group; Section~\ref{sec:proof-main} completes the proof.

\section{Enriques lattices and metric periods}
\label{sec:metric-periods}

\subsection{The K3 cover and the Enriques lattices}
\label{subsec:enriques-lattices}

Fix the canonical cover $p:X\to E$ with covering involution $\iota$, and give both manifolds their complex orientations. Write
$$L:=U\oplus E_8(-1),\qquad\Lambda:=U\oplus L\oplus L,\qquad\tau(u,x,y):=(-u,y,x).$$
Here $U$ is the hyperbolic plane, $E_8$ is positive definite, and $K(k)$ means that the bilinear form of $K$ is multiplied by $k$. We use $K_{\mathbb{F}}:=K\otimes_{\mathbb{Z}}\mathbb{F}$. Fix a compatible marking
$$\alpha:H^2(X;\mathbb{Z})\xrightarrow{\sim}\Lambda,\qquad\alpha\circ\iota^*=\tau\circ\alpha;$$
see \cite[Section~2.1]{PetersSterk2020}. The eigensublattices are primitive and orthogonal, with
$$M=\{(0,x,x)\mid x\in L\}\cong L(2),$$
$$N=\{(u,x,-x)\mid u\in U,\ x\in L\}\cong U\oplus L(2).$$
They have signatures $(1,9)$ and $(2,10)$, and $N=M^\perp$. Pullback identifies $H^2(E;\mathbb{Z})/\operatorname{Tor}$ with $M$, multiplying the intersection pairing by two. Although $\Lambda_{\mathbb{R}}=M_{\mathbb{R}}\oplus N_{\mathbb{R}}$, the integral sum has index $2^{10}$.

For an even lattice $K$, let $K^\vee$ be its dual, $A_K:=K^\vee/K$, and $q_K(x+K):=(x,x)\bmod2\mathbb{Z}$. Unimodularity of $\Lambda$ and primitivity of $M,N$ identify $\Lambda/(M\oplus N)$ with the graph of an anti-isometry $\varphi:A_M\to A_N$. Explicitly,
$$A_M\cong A_N\cong(\mathbb{Z}/2\mathbb{Z})^{10},\qquad\Lambda=\{(m,n)\in M^\vee\oplus N^\vee\mid n+N=\varphi(m+M)\}.$$
For $g_K\in\operatorname{O}(K)$, denote its discriminant action by $\overline g_K$.

\begin{lemma}[Extension criterion]
\label{lem:gluing-extension}
A pair $(g_M,g_N)\in\operatorname{O}(M)\times\operatorname{O}(N)$ extends to an isometry of $\Lambda$ if and only if
$$\overline g_N\circ\varphi=\varphi\circ\overline g_M.$$
The extension is unique and commutes with $\tau$.
\end{lemma}

\begin{proof}
The rational isometry $g_M\oplus g_N$ preserves $\Lambda$ precisely when it preserves the graph of $\varphi$. Uniqueness follows from the finite index of $M\oplus N$; preservation of the eigenspaces gives commutation with $\tau$.
\end{proof}

\subsection{Positive planes and period domains}
\label{subsec:period-domains}

Let $\operatorname{Gr}_r^+(V)$ denote the space of positive $r$-planes in a real quadratic space $V$, without an orientation. Define
$$\mathcal{D}_M:=\operatorname{Gr}_1^+(M_{\mathbb{R}}),\qquad\mathcal{D}_N:=\operatorname{Gr}_2^{+,\mathrm{or}}(N_{\mathbb{R}}).$$
The second space is the type IV domain
$$\mathcal{D}_N\cong\{[\omega]\in\mathbb{P}(N_{\mathbb{C}})\mid(\omega,\omega)=0,\ (\omega,\overline\omega)>0\}.$$
The period line determines the plane spanned and oriented by $(\operatorname{Re}\omega,\operatorname{Im}\omega)$. Orientation reversal is $c([\omega])=[\overline\omega]$. It exchanges the two components of $\mathcal{D}_N$; fix one, denoted $\mathcal{D}_N^\circ$, and write $\underline P$ for the unoriented plane underlying $P$.

\begin{lemma}
\label{lem:period-domain-topology}
The spaces $\mathcal{D}_M$ and $\mathcal{D}_N^\circ$ are contractible, of real dimensions $9$ and $20$. Forgetting orientation identifies $\mathcal{D}_N^\circ$ diffeomorphically with $\operatorname{Gr}_2^+(N_{\mathbb{R}})$.
\end{lemma}

\begin{proof}
For a quadratic space of signature $(r,s)$, fix an orthogonal positive-negative decomposition $V=V_+\oplus V_-$. Every maximal positive plane is the graph of a unique map $T:V_+\to V_-$ with $\|T\|_{\mathrm{op}}<1$. This is a convex open set of dimension $rs$. The two orientations of $V_+$ give the two components of the oriented space.
\end{proof}

Set $\mathcal{D}_\Lambda:=\operatorname{Gr}_3^+(\Lambda_{\mathbb{R}})$ and let $\mathcal{D}_\Lambda^\tau$ be its locus of setwise $\tau$-invariant planes.

\begin{lemma}
\label{lem:invariant-three-planes}
The map
$$\Psi:\mathcal{D}_M\times\mathcal{D}_N^\circ\xrightarrow{\sim}\mathcal{D}_\Lambda^\tau,\qquad(\ell,P)\longmapsto\ell\oplus\underline P$$
is a diffeomorphism.
\end{lemma}

\begin{proof}
An invariant positive three-plane splits into its intersections with $M_{\mathbb{R}}$ and $N_{\mathbb{R}}$. Their positive indices force the dimensions to be one and two. The choice of $\mathcal{D}_N^\circ$ supplies the unique required orientation; orthogonal sum is the inverse construction.
\end{proof}

The choice of sheet is only an orientation convention for metric periods. If $g_N\in\operatorname{O}(N)$ exchanges the sheets, its action on the chosen sheet is $P\mapsto c(g_NP)$, an anti-holomorphic map.

\subsection{Metric periods and discriminants}
\label{subsec:metric-discriminants}

Fix a connected component $\mathcal{T}_E$ of
$$\mathcal{T}^{\mathrm{Ric}=0}_{\mathrm{unit}}(E):=\mathcal{R}^{\mathrm{Ric}=0}_{\mathrm{unit}}(E)/\operatorname{Diff}_0(E).$$
For a unit-volume metric $g$ on $E$, the normalized lift $\widehat g:=2^{-1/2}p^*g$ has volume one. Its self-dual harmonic classes span a $\tau$-invariant positive three-plane $H_g$. Thus $\Psi^{-1}(H_g)$ defines a period map $\operatorname{Per}_E:\mathcal{T}_E\to\mathcal{D}_M\times\mathcal{D}_N^\circ$.

Put $\operatorname{Roots}(K):=\{\delta\in K\mid\delta^2=-2\}$. For $\delta\in\operatorname{Roots}(\Lambda)$, its rational projections are $\delta_M=(\delta+\tau\delta)/2$ and $\delta_N=(\delta-\tau\delta)/2$, and its wall is
$$Z_\delta:=\{(\ell,P)\mid\ell\perp\delta_M,\ \underline P\perp\delta_N\}.$$
For $\delta\in N$, write $D_\delta:=\{P\in\mathcal{D}_N^\circ\mid\underline P\perp\delta\}$, so $Z_\delta=\mathcal{D}_M\times D_\delta$. Define
$$\Delta_{\mathrm{per}}:=\bigcup_{\delta\in\operatorname{Roots}(N)}\mathcal{D}_M\times D_\delta,\qquad\Delta_{\mathrm{met}}:=\bigcup_{\delta\in\operatorname{Roots}(\Lambda)\setminus N}Z_\delta,$$
$$\mathcal{U}_E:=(\mathcal{D}_M\times\mathcal{D}_N^\circ)\setminus(\Delta_{\mathrm{per}}\cup\Delta_{\mathrm{met}}).$$

\begin{lemma}
\label{lem:discriminant-codimensions}
The nonempty walls form a locally finite family of closed smooth submanifolds. Their real codimension is two for roots in $N$ and three for the remaining roots. The complement $\mathcal{U}_E$ is path connected.
\end{lemma}

\begin{proof}
For a root in $N$, its orthogonal complement has signature $(2,9)$, so $D_\delta$ has dimension $18$. Otherwise both projections of $\delta$ are nonzero: the $N$-projection cannot vanish because $M=L(2)$ has no roots. If $Z_\delta$ is nonempty, both projections have negative square. The two orthogonal complements then have signatures $(1,8)$ and $(2,9)$, giving a wall of dimension $8+18=26$ in a domain of dimension $29$.

For local finiteness, a positive three-plane $H$ defines the positive-definite form
$$Q_H(v):=(v_H,v_H)-(v_{H^\perp},v_{H^\perp}).$$
A root perpendicular to $H$ has $Q_H(\delta)=2$. On a compact family of planes, these forms uniformly bound a fixed Euclidean norm from below. Only finitely many lattice vectors can therefore define a wall meeting that family. Finally, relative transversality perturbs paths off this locally finite family of submanifolds of codimension at least two.
\end{proof}

\begin{proposition}[Metric period description]
\label{prop:metric-period-description}
The period map is a homeomorphism $\operatorname{Per}_E:\mathcal{T}_E\xrightarrow{\sim}\mathcal{U}_E$.
\end{proposition}

\begin{proof}
Choose $[g_0]\in\mathcal{T}_E$ and let $\mathcal{C}_X$ be the K3 metric Teichm\"uller component containing $[\widehat g_0]$. Raman's lifting results identify $\mathcal{T}_E$ with the component of $\mathcal{C}_X^{[\iota]}$ containing this point; see \cite[Corollaries~6.11--6.12 and the proof of Corollary~6.12]{Raman2025}.

The K3 metric Torelli theorem \cite[Proposition~2.1]{FarbLooijenga2024} identifies $\mathcal{C}_X$ with
$$\mathcal{U}_\Lambda:=\{H\in\mathcal{D}_\Lambda\mid H\not\subset\delta^\perp\text{ for every }\delta\in\operatorname{Roots}(\Lambda)\}.$$
Equivariance identifies the fixed locus with $\mathcal{U}_\Lambda\cap\mathcal{D}_\Lambda^\tau$, whose inverse image under $\Psi$ is $\mathcal{U}_E$. This space is connected by Lemma~\ref{lem:discriminant-codimensions}, so the entire fixed locus is the required component.
\end{proof}

\subsection{The diffeomorphism-induced action}
\label{subsec:actual-action}

Use the left action $f\cdot g:=(f^{-1})^*g$, and put $\operatorname{Mod}(E):=\operatorname{Diff}(E)/\operatorname{Diff}_0(E)$. Let $\operatorname{Mod}(E)_{\mathcal{T}_E}$ be the subgroup preserving $\mathcal{T}_E$. The ambient gluing-compatible group is
$$\Gamma_E^{\mathrm{gl}}:=\{(g_M,g_N)\in\operatorname{O}(M)\times\operatorname{O}(N)\mid\overline g_N\circ\varphi=\varphi\circ\overline g_M\}.$$
Let $\varepsilon_N(g_N)\in\mathbb{Z}/2\mathbb{Z}$ record whether $g_N$ exchanges the sheets. The action on our product domain is
$$(g_M,g_N)\cdot(\ell,P):=(g_M\ell,c^{\varepsilon_N(g_N)}g_NP).$$
It preserves both discriminants: the extension $A$ to $\Lambda$ sends $Z_\delta$ to $Z_{A\delta}$.

Every self-diffeomorphism of $E$ preserves orientation, since its signature is nonzero, and lifts to $X$. The two lifts commute with $\iota$ and differ by $\iota$. A lift $\widetilde f$ gives
$$A_f:=\alpha\circ(\widetilde f^{-1})^*\circ\alpha^{-1}.$$
Changing the lift multiplies $A_f$ by $\tau$, which acts trivially on metric periods. Let $\widetilde\Gamma_E^{\mathrm{act}}\subset\Gamma_E^{\mathrm{gl}}$ consist of these restrictions for $[f]\in\operatorname{Mod}(E)_{\mathcal{T}_E}$, allowing both lifts, and define
$$\Gamma_E^{\mathrm{act}}:=\operatorname{im}\bigl(\operatorname{Mod}(E)_{\mathcal{T}_E}\longrightarrow\operatorname{Homeo}(\mathcal{U}_E)\bigr).$$
Naturality of harmonic forms identifies this with the effective image of $\widetilde\Gamma_E^{\mathrm{act}}$. Actual lifts lie in $\operatorname{O}^+(\Lambda)$, the subgroup preserving the orientation of maximal positive three-planes \cite{FarbLooijenga2024}.

\begin{lemma}
\label{lem:actual-action-kernel}
There is an exact sequence
$$1\longrightarrow\langle\tau\rangle\longrightarrow\widetilde\Gamma_E^{\mathrm{act}}\longrightarrow\Gamma_E^{\mathrm{act}}\longrightarrow1.$$
\end{lemma}

\begin{proof}
If an actual pair acts trivially on $\mathcal{U}_E$, density and continuity show that its restrictions fix all positive lines in $M_{\mathbb{R}}$ and all positive two-planes in $N_{\mathbb{R}}$. These restrictions are scalar: a positive line in $N_{\mathbb{R}}$ is the intersection of the positive two-planes containing it. Thus $g_M=\pm\operatorname{id}_M$ and $g_N=\pm\operatorname{id}_N$. The negative sign on $M$ reverses the orientation of a positive three-plane and is excluded. The remaining pairs are $1$ and $\tau$, both of which act trivially.
\end{proof}

The sheet character descends to $\varepsilon_E:\Gamma_E^{\mathrm{act}}\to\mathbb{Z}/2\mathbb{Z}$. Write $\Gamma_E^+:=\ker\varepsilon_E$ for the sheet-preserving subgroup.

\begin{proposition}[The coarse quotient and its stabilizers]
\label{prop:coarse-metric-quotient}
The coarse component $\mathcal{M}_0$ containing the image of $\mathcal{T}_E$ is naturally homeomorphic to $\Gamma_E^{\mathrm{act}}\backslash\mathcal{U}_E$. The action is properly discontinuous with finite stabilizers, and
$$\operatorname{Isom}(E,g)\xrightarrow{\sim}(\Gamma_E^{\mathrm{act}})_{\operatorname{Per}_E([g])}.$$
\end{proposition}

\begin{proof}
The quotient by $\operatorname{Diff}(E)$ is obtained by first dividing by $\operatorname{Diff}_0(E)$ and then by $\operatorname{Mod}(E)$. Components of metric Teichm\"uller space are open by local deformation theory \cite[Section~6]{Raman2025}. The union of the translates of $\mathcal{T}_E$ is therefore an open and closed saturated subset, and its connected quotient is exactly $\mathcal{M}_0$. Orbit maps are open, so the period homeomorphism gives the asserted homeomorphism of quotients.

The product of the maximal positive-plane domains is homogeneous with compact isotropy. Its real orthogonal group acts properly; the integral lifted group is discrete. Proper discontinuity and finiteness of stabilizers pass to its effective image and the invariant open subset $\mathcal{U}_E$.

If $f$ represents a period stabilizer, injectivity of $\operatorname{Per}_E$ gives $(f^{-1})^*g=u^*g$ for some $u\in\operatorname{Diff}_0(E)$. Hence $u\circ f$ is an isometry inducing the same period action. Conversely, an isometry acting trivially on periods has a lift acting trivially on $H^2(X;\mathbb{Z})$, after composing with $\iota$ if necessary. Such a K3 isometry is the identity \cite[Corollary~2.2]{FarbLooijenga2024}, proving injectivity.
\end{proof}

\section{Discriminant complements and action braid groups}
\label{sec:action-braid-groups}

\subsection{Meridians in the period complement}
\label{subsec:period-complement}

Put
$$\mathcal{U}_N:=\mathcal{D}_N^\circ\setminus\bigcup_{\delta\in\operatorname{Roots}(N)}D_\delta,\qquad\mathcal{V}_E:=\mathcal{D}_M\times\mathcal{U}_N.$$
Thus $\mathcal{U}_E=\mathcal{V}_E\setminus\Delta_{\mathrm{met}}$. Fix $x_0=(\ell_0,P_0)\in\mathcal{U}_E$ and $t_0:=\operatorname{Per}_E^{-1}(x_0)$.

\begin{proposition}[Removing the metric discriminant]
\label{prop:codimension-three-removal}
Inclusion and projection induce isomorphisms
$$\pi_1(\mathcal{T}_E,t_0)\cong\pi_1(\mathcal{U}_E,x_0)\xrightarrow{\sim}\pi_1(\mathcal{V}_E,x_0)\xrightarrow{\sim}\pi_1(\mathcal{U}_N,P_0).$$
\end{proposition}

\begin{proof}
The removed walls are locally finite closed smooth submanifolds of real codimension three. Relative transversality perturbs paths and disk homotopies off these walls while fixing their endpoints or boundary. Compactness reduces each perturbation to finitely many walls, proving both surjectivity and injectivity. Projection is a homotopy equivalence because $\mathcal{D}_M$ is contractible.
\end{proof}

Each root mirror $D_\delta$ is connected. Distinct mirrors meet, if at all, in complex codimension two in $\mathcal{D}_N^\circ$: their two roots span a negative-definite plane. Removing these intersections from a mirror leaves a nonempty path-connected regular part. A small transverse complex disk at a regular point defines a positive meridian $\mu_\delta$ based at $P_0$, well defined up to conjugation. The roots $\delta$ and $-\delta$ give the same class.

\begin{proposition}[Meridians and first homology]
\label{prop:meridian-homology}
The meridians normally generate $\pi_1(\mathcal{T}_E,t_0)$, and
$$H_1(\mathcal{T}_E;\mathbb{Z})\cong\bigoplus_{[\delta]\in\operatorname{Roots}(N)/\{\pm1\}}\mathbb{Z}[\mu_\delta].$$
This group has countably infinite rank. In particular, $\mathcal{T}_E$ is not simply connected.
\end{proposition}

\begin{proof}
Work in $\mathcal{U}_N$ using Proposition~\ref{prop:codimension-three-removal}. A loop bounds a disk in the contractible domain $\mathcal{D}_N^\circ$. Perturb the disk relative to its boundary to meet the mirrors transversely and avoid their pairwise intersections. Its finitely many intersection points give a decomposition of the loop into conjugates of meridians.

For independence in homology, trivialize the tautological line bundle smoothly over $\mathcal{D}_N^\circ$ by a nowhere vanishing section $s$. The function $f_\delta(P):=(s(P),\delta)$ has zero locus $D_\delta$, with multiplicity one in a local holomorphic trivialization. On $\mathcal{U}_N$ its winding number is one on $\mu_\delta$ and zero on every other mirror meridian. These homomorphisms detect every finite relation among the meridian classes and prove the direct-sum formula.

There are countably infinitely many mirrors. Indeed, for a hyperbolic basis $e,f$ of $U$ and a nonzero isotropic vector $u\in U(2)$, the roots $e-f+ku$, $k\in\mathbb{Z}$, are distinct modulo sign.
\end{proof}

\subsection{Action braid groups and the coarse fundamental group}
\label{subsec:action-braid-groups}

Let a discrete group $\Gamma$ act effectively and properly discontinuously on a connected smooth manifold $Z$, and fix $z_0\in Z$. Proper discontinuity here means that every compact set meets only finitely many of its translates. Define
$$\Pi_1^\Gamma(Z,z_0):=\{[\alpha,\gamma]\mid\gamma\in\Gamma,\ \alpha:z_0\to\gamma z_0\},\qquad[\alpha,\gamma][\beta,\eta]:=[\alpha*(\gamma\beta),\gamma\eta].$$
Paths are taken up to homotopy relative to endpoints, and $*$ means traversal from left to right. The endpoint map gives an exact sequence
$$1\longrightarrow\pi_1(Z,z_0)\longrightarrow\Pi_1^\Gamma(Z,z_0)\xrightarrow{e_Z}\Gamma\longrightarrow1.$$
For $\gamma z=z$ and $\beta:z_0\to z$, the stabilizer class is
$$s(\gamma,z,\beta):=[\beta*(\gamma\beta)^{-1},\gamma].$$
For fixed $z,\beta$, this embeds $\Gamma_z$ as a subgroup. Let $\mathcal{N}_Z^{\mathrm{stab}}$ be the normal subgroup generated by all such classes.

\begin{proposition}[Armstrong's theorem in action-group form]
\label{prop:armstrong-action-group}
Projection of paths to the coarse quotient induces
$$\pi_1(\Gamma\backslash Z,[z_0])\cong\Pi_1^\Gamma(Z,z_0)/\mathcal{N}_Z^{\mathrm{stab}}.$$
\end{proposition}

\begin{proof}
In the path model of the universal cover $\widetilde Z$, the action
$$[\alpha,\gamma]\cdot[\lambda]:=[\alpha*(\gamma\lambda)]$$
identifies $\Pi_1^\Gamma(Z,z_0)$ with the group $\widehat\Gamma$ of all lifts of $\Gamma$ to $\widetilde Z$. This action is properly discontinuous: for a compact set upstairs, proper discontinuity downstairs leaves finitely many possible endpoint elements, and proper discontinuity of the deck group leaves finitely many lifts of each. Moreover, $\widehat\Gamma\backslash\widetilde Z\cong\Gamma\backslash Z$.

Armstrong's theorem \cite{Armstrong1968} identifies the fundamental group with $\widehat\Gamma$ modulo the normal closure of its fixed-point elements. A lift $[\alpha,\gamma]$ fixes $[\beta]$ precisely when $\gamma\beta(1)=\beta(1)$ and $\alpha\simeq\beta*(\gamma\beta)^{-1}$ relative to endpoints. These are exactly the stated stabilizer classes.
\end{proof}

Using the period homeomorphism, put
$$\mathcal{B}_E:=\Pi_1^{\Gamma_E^{\mathrm{act}}}(\mathcal{U}_E,x_0)\cong\Pi_1^{\Gamma_E^{\mathrm{act}}}(\mathcal{T}_E,t_0).$$
Write $e_E$ for its endpoint map, and $\mathcal{N}_{\mathrm{stab}}\triangleleft\mathcal{B}_E$ for the normal closure of stabilizer classes at points of $\mathcal{U}_E$.

\begin{corollary}[The coarse metric fundamental group]
\label{cor:coarse-braid-quotient}
For the corresponding coarse component,
$$\pi_1(\mathcal{M}_0)\cong\mathcal{B}_E/\mathcal{N}_{\mathrm{stab}}.$$
\end{corollary}

\begin{proof}
Apply Proposition~\ref{prop:armstrong-action-group} to the action in Proposition~\ref{prop:coarse-metric-quotient}.
\end{proof}

Similarly, define $\mathcal{B}_E^+:=e_E^{-1}(\Gamma_E^+)$ and let $\mathcal{N}_{\mathrm{stab}}^+$ be the normal closure, in this subgroup, of sheet-preserving stabilizer classes. Then
$$\pi_1(\Gamma_E^+\backslash\mathcal{T}_E)\cong\mathcal{B}_E^+/\mathcal{N}_{\mathrm{stab}}^+.$$
By Proposition~\ref{prop:coarse-metric-quotient}, all these stabilizers are isometries of smooth metrics.

\begin{lemma}[Equivariant deletion]
\label{lem:equivariant-deletion}
Inclusion induces an isomorphism
$$\Pi_1^{\Gamma_E^{\mathrm{act}}}(\mathcal{U}_E,x_0)\xrightarrow{\sim}\Pi_1^{\Gamma_E^{\mathrm{act}}}(\mathcal{V}_E,x_0),$$
and likewise for $\Gamma_E^+$.
\end{lemma}

\begin{proof}
For each endpoint element $\gamma$, the relative transversality argument of Proposition~\ref{prop:codimension-three-removal} identifies the path classes from $x_0$ to $\gamma x_0$ in the two spaces. Equivariance makes these identifications compatible with multiplication.
\end{proof}

This last isomorphism concerns action braid groups. Stabilizer relations must still come from $\mathcal{U}_E$: a fixed point in $\mathcal{V}_E$ may lie on the metric discriminant.

\section{Smooth metric stabilizers}
\label{sec:metric-stabilizers}

\subsection{Realization criteria}
\label{subsec:realization-criteria}

Let $K\subset\operatorname{O}^+(\Lambda)$ be finite, contain $\tau$, and centralize $\tau$. Its invariant-plane domain is
$$\mathcal{D}_\Lambda^K:=\{H\in\mathcal{D}_\Lambda\mid aH=H\text{ for all }a\in K\}\subset\mathcal{D}_\Lambda^\tau.$$
Invariance is setwise.

\begin{lemma}[Realization at a smooth period]
\label{lem:realization-at-period}
For $H\in\mathcal{D}_\Lambda^K\cap\mathcal{U}_\Lambda$, let $x:=\Psi^{-1}(H)$ and choose a metric representing $\operatorname{Per}_E^{-1}(x)$. There is an injection
$$K/\langle\tau\rangle\longrightarrow\operatorname{Isom}(E,g)$$
whose period action is the prescribed lattice action. In particular, it belongs to $\Gamma_E^{\mathrm{act}}$.
\end{lemma}

\begin{proof}
By \cite[Corollary~2.2]{FarbLooijenga2024}, each $a\in K$ is represented by a unique isometry $h_a$ of $(X,\widehat g)$ with $\alpha\circ(h_a^{-1})^*\circ\alpha^{-1}=a$. Uniqueness makes this an injective group homomorphism. The isometry corresponding to $\tau$ is the original deck involution $\iota$, so all $h_a$ commute with $\iota$. Descent has kernel exactly $\langle\iota\rangle$. The descended isometries fix $[g]$, and hence preserve its component $\mathcal{T}_E$.
\end{proof}

A $K$-invariant positive three-plane $H_0$ exists: average a Euclidean inner product over $K$ and take the positive spectral subspace of the self-adjoint operator representing the lattice form. Graph coordinates then give
$$\mathcal{D}_\Lambda^K\cong\{T\in\operatorname{Hom}_K(H_0,H_0^\perp)\mid\|T\|_{\mathrm{op}}<1\}.$$
Define
$$W_K:=\operatorname{span}_{\mathbb{R}}\{H\mid H\in\mathcal{D}_\Lambda^K\},\qquad L_K:=\Lambda\cap W_K^\perp.$$
The lattice $L_K$ is negative definite. Moreover, $W_K=H_0\oplus\sum_T\operatorname{im}T$ is the sum of the real isotypical summands whose irreducible types occur in $H_0$. Thus this is the realization lattice in the metric criterion of \cite[Theorem~1.2]{FarbLooijenga2024}.

\begin{proposition}[Root avoidance and metric realization]
\label{prop:finite-group-realization}
The intersection $\mathcal{D}_\Lambda^K\cap\mathcal{U}_\Lambda$ is nonempty if and only if $L_K$ contains no roots. In that case it is open and dense, and every plane in it realizes $K/\langle\tau\rangle$ as isometries of a smooth metric in the prescribed component $\mathcal{T}_E$.
\end{proposition}

\begin{proof}
A root in $L_K$ is perpendicular to every invariant plane, proving necessity. Conversely, a root wall in the graph coordinates is given by
$$(\delta,v+Tv)=0\qquad\text{for all }v\in H_0.$$
If $L_K$ is root-free, these affine linear equations do not hold throughout the domain. Each wall therefore has empty interior. Local finiteness of the root walls gives an open dense complement. Proposition~\ref{prop:metric-period-description} and Lemma~\ref{lem:realization-at-period} realize every point of this complement in $\mathcal{T}_E$.
\end{proof}

\begin{lemma}[Averaging a K\"ahler class]
\label{lem:averaging-kahler-class}
A finite group $F$ of holomorphic automorphisms of an Enriques surface preserves a smooth unit-volume Ricci-flat K\"ahler metric.
\end{lemma}

\begin{proof}
For a K\"ahler class $\kappa_0$, the average $|F|^{-1}\sum_{f\in F}f^*\kappa_0$ is K\"ahler and $F$-invariant. The uniqueness in Yau's theorem \cite{Yau1978} makes its Ricci-flat representative invariant. Rescale to volume one.
\end{proof}

\begin{corollary}[Transfer of finite metric actions]
\label{cor:transfer-finite-actions}
Let $F\subset\operatorname{Isom}(E,g_0)$ be finite, and let $K$ be the marked cohomological image of its full group of lifts to $X$. There is an isomorphic group $F'\subset\operatorname{Isom}(E,g)$ for a metric $[g]\in\mathcal{T}_E$, at the same marked period, whose full lifted cohomological image is $K$.
\end{corollary}

\begin{proof}
The lifts act faithfully on K3 cohomology and contain the deck involution, so $F\cong K/\langle\tau\rangle$. Their common metric period $H_{g_0}$ is root-free. The period homeomorphism on $\mathcal{T}_E$ and Lemma~\ref{lem:realization-at-period} realize the same $K$ there. Since the realized group contains the deck involution, it contains both lifts of every descended element.
\end{proof}

\subsection{Bertini stabilizers and the projective kernel}
\label{subsec:bertini-kernel}

Choose a component of the positive cone of $L_{\mathbb{R}}$, with stabilizer $\operatorname{O}^{\uparrow}(L)$. On $V:=L/2L$, put $q_2(\bar v):=(v,v)_L/2\bmod2$, and define
$$W(2):=\ker\bigl(\operatorname{O}^{\uparrow}(L)\longrightarrow\operatorname{O}(V,q_2)\bigr).$$
The identification $V\to A_M$, $\bar v\mapsto v/2+M$, respects orthogonal-group actions. Thus $W(2)$ consists of the cone-preserving isometries of $M$ acting trivially on $A_M$.

For a decomposition $L=U'\oplus E'$ with $U'\cong U$ and $E'\cong E_8(-1)$, let $b$ be the identity on $U'$ and minus the identity on $E'$. Such Bertini involutions generate $W(2)$ by \cite[Theorem~1.1 and its proof]{Allcock2018}.

\begin{lemma}[Bertini metric stabilizers]
\label{lem:bertini-metric-stabilizers}
Every pair $(b,\operatorname{id}_N)$ extends to an involution $\widetilde b\in\operatorname{O}^+(\Lambda)$ commuting with $\tau$. Its effective action $\gamma_b$ is a sheet-preserving stabilizer of a smooth metric in $\mathcal{T}_E$.
\end{lemma}

\begin{proof}
The discriminant action of $b$ is trivial, so the extension criterion applies. The extension fixes a positive three-plane pointwise and hence lies in $\operatorname{O}^+(\Lambda)$. For $K_b:=\langle\tau,\widetilde b\rangle$, invariant planes consist of a positive line in $U'(2)_{\mathbb{R}}$ and an arbitrary positive two-plane in $N_{\mathbb{R}}$. Consequently,
$$W_{K_b}=U'(2)_{\mathbb{R}}\oplus N_{\mathbb{R}},\qquad L_{K_b}=\Lambda\cap E'(2)_{\mathbb{R}}=E'(2)\cong E_8(-2).$$
The integral equality follows from $\Lambda\cap M_{\mathbb{R}}=M$ and the direct-sum decomposition of $M$. This lattice has no roots, so Proposition~\ref{prop:finite-group-realization} gives the required smooth metric.
\end{proof}

Let $\operatorname{O}^+(N)$ preserve $\mathcal{D}_N^\circ$, and set $Q_N:=\operatorname{O}^+(N)/\{\pm\operatorname{id}_N\}$. Restriction defines
$$\rho_N:\Gamma_E^+\longrightarrow Q_N,\qquad\gamma\longmapsto[g_N],\qquad Q_N^{\mathrm{act}}:=\rho_N(\Gamma_E^+).$$
The two lifts differ by $\tau$, so the projective class is well defined.

\begin{proposition}[Kernel of projective restriction]
\label{prop:projective-period-kernel}
There is an exact sequence
$$1\longrightarrow W(2)\longrightarrow\Gamma_E^+\xrightarrow{\rho_N}Q_N^{\mathrm{act}}\longrightarrow1.$$
The first map sends $w$ to the transformation represented by $(w,\operatorname{id}_N)$, and its image is generated by the $\gamma_b$.
\end{proposition}

\begin{proof}
A lift of an element of $\ker\rho_N$ can be uniquely chosen with $g_N=\operatorname{id}_N$. Preservation of positive-three-plane orientation then forces $g_M$ to preserve the positive cone component. Gluing gives trivial action on $A_M$, so $g_M\in W(2)$. Conversely, every $w\in W(2)$ is a product of Bertini involutions, each realized in the chosen component by Lemma~\ref{lem:bertini-metric-stabilizers}. The action is faithful by Lemma~\ref{lem:actual-action-kernel}.
\end{proof}

Projection $p_N:\mathcal{U}_E\to\mathcal{U}_N$ is equivariant for $\rho_N$ and induces an isomorphism on ordinary fundamental groups. Set $\mathcal{B}_N^{\mathrm{act}}:=\Pi_1^{Q_N^{\mathrm{act}}}(\mathcal{U}_N,P_0)$ and define
$$p_{\#}:\mathcal{B}_E^+\longrightarrow\mathcal{B}_N^{\mathrm{act}},\qquad[\alpha,\gamma]\longmapsto[p_N\circ\alpha,\rho_N(\gamma)].$$
Let $\mathcal{N}_{\mathrm{Bert}}\triangleleft\mathcal{B}_E^+$ be normally generated by all Bertini stabilizer classes $s(\gamma_b,x,\beta)$ at smooth metric periods.

\begin{theorem}[Exact removal of the Bertini kernel]
\label{thm:bertini-kernel-removal}
The map $p_{\#}$ is surjective, and
$$\ker p_{\#}=\mathcal{N}_{\mathrm{Bert}},\qquad\mathcal{B}_E^+/\mathcal{N}_{\mathrm{Bert}}\cong\mathcal{B}_N^{\mathrm{act}}.$$
Moreover, the endpoint map restricts to an isomorphism $\ker p_{\#}\xrightarrow{\sim}\ker\rho_N\cong W(2)$.
\end{theorem}

\begin{proof}
Given a period braid, lift its endpoint through $\rho_N$ and choose any path to that endpoint in $\mathcal{U}_E$. The discrepancy between its projection and the prescribed period path is a loop. Surjectivity of $(p_N)_*$ allows us to correct the chosen path by a based loop, proving surjectivity of $p_{\#}$. The same correction realizes every element of $\ker\rho_N$ as the endpoint of an element of $\ker p_{\#}$. Injectivity on this kernel follows from injectivity of $(p_N)_*$: a kernel element with trivial endpoint is a loop whose projection is null-homotopic.

For a Bertini stabilizer, $p_N\circ\gamma_b=p_N$, so
$$p_{\#}s(\gamma_b,x,\beta)=[(p_N\circ\beta)*(p_N\circ\beta)^{-1},1]=1.$$
Under the endpoint isomorphism just proved, one such stabilizer for each $b$ maps to a generating set of $W(2)$. Thus these classes generate $\ker p_{\#}$, proving equality with their normal closure.
\end{proof}

\subsection{The full projective period action}
\label{subsec:projective-period-action}

Write $\operatorname{O}^{\uparrow}(M)$ for the positive-cone-component stabilizer and put
$$\operatorname{O}^{+,+}(\Lambda)^\tau:=\{a\in\operatorname{O}(\Lambda)\mid a\tau=\tau a,\ a|_M\in\operatorname{O}^{\uparrow}(M),\ a|_N\in\operatorname{O}^+(N)\}.$$
This group lies in $\operatorname{O}^+(\Lambda)$ and contains $\langle\tau\rangle$. Choose $[g_0]\in\mathcal{T}_E$ and an Enriques complex structure $J_0$ for which $g_0$ is K\"ahler. Oberdieck identifies the marked lifted monodromy group of Enriques surfaces, modulo the deck involution, as
$$\widetilde{\operatorname{Mon}}(E,J_0)=\operatorname{O}^{+,+}(\Lambda)^\tau/\langle\tau\rangle;$$
see \cite[Proposition~3.3]{Oberdieck2025}.

\begin{lemma}[Monodromy preserves the metric component]
\label{lem:monodromy-metric-component}
Each lifted monodromy class has a diffeomorphism representative preserving $\mathcal{T}_E$, whose period action belongs to $\Gamma_E^+$.
\end{lemma}

\begin{proof}
For a monodromy generator, take a smooth projective family along a path, with both endpoint fibres identified biholomorphically with $(E,J_0)$. A smooth trivialization $F_t$ along the path starts at the first identification and ends at the second composed with a diffeomorphism $h:E\to E$. Parallel transport on the K3 covers is represented by $(\widetilde h^{-1})^*$, up to the deck involution.

A relative polarization gives fibrewise K\"ahler classes. Their normalized Ricci-flat representatives $r_t$ vary smoothly: after normalizing the potential, the linearization of the parameter-dependent Monge--Amp\`ere equation is the Laplacian on mean-zero functions, so the implicit function theorem, elliptic regularity, and uniqueness apply. Pulling back gives a metric path from a $J_0$-K\"ahler metric $g_a$ to $h^*g_b$, where $g_b$ is also $J_0$-K\"ahler.

Convexity of the K\"ahler cone and the same dependence on the class give paths from $g_0$ to $g_a$ and from $g_b$ to $g_0$. Concatenating with the pulled-back family produces a path from $g_0$ to $h^*g_0$. Thus both $h$ and $h^{-1}$ preserve $\mathcal{T}_E$. The monodromy theorem ensures sheet preservation, and products of generators retain these properties.
\end{proof}

\begin{proposition}[The sheet-preserving arithmetic action]
\label{prop:sheet-preserving-arithmetic-action}
The lattice action induces an isomorphism
$$\operatorname{O}^{+,+}(\Lambda)^\tau/\langle\tau\rangle\xrightarrow{\sim}\Gamma_E^+.$$
Furthermore, $Q_N^{\mathrm{act}}=Q_N$, so
$$1\longrightarrow W(2)\longrightarrow\Gamma_E^+\xrightarrow{\rho_N}Q_N\longrightarrow1.$$
\end{proposition}

\begin{proof}
An actual sheet-preserving lift preserves positive-three-plane orientation, hence also the positive cone component in $M$. It belongs to $\operatorname{O}^{+,+}(\Lambda)^\tau$. Conversely, Oberdieck's theorem and Lemma~\ref{lem:monodromy-metric-component} give a representative of every class preserving the prescribed component. Lemma~\ref{lem:actual-action-kernel} identifies the kernel.

For $u\in\operatorname{O}^+(N)$, surjectivity of reduction
$$\operatorname{O}^{\uparrow}(L)\longrightarrow\operatorname{O}(V,q_2)$$
from \cite[Sections~1--2]{Allcock2018} provides $v\in\operatorname{O}^{\uparrow}(M)$ with $\bar v=\varphi^{-1}\bar u\varphi$. The pair $(v,u)$ extends to $\operatorname{O}^{+,+}(\Lambda)^\tau$, and the first assertion makes its action actual. Its image is $[u]$, proving surjectivity; the kernel was computed in Proposition~\ref{prop:projective-period-kernel}.
\end{proof}

Henceforth write $\mathcal{B}_N:=\Pi_1^{Q_N}(\mathcal{U}_N,P_0)=\mathcal{B}_N^{\mathrm{act}}$. The Bertini reduction becomes
$$1\longrightarrow\mathcal{N}_{\mathrm{Bert}}\longrightarrow\mathcal{B}_E^+\xrightarrow{p_{\#}}\mathcal{B}_N\longrightarrow1.$$
This identifies the effective action on the chosen component; no triviality of the smooth Torelli group is asserted.

\subsection{Kantor stabilizers and retained reflection mirrors}
\label{subsec:kantor-stabilizers}

For a primitive vector $n\in N$, let $\operatorname{div}_N(n)$ be the positive generator of the ideal $(n,N)\subset\mathbb{Z}$. Put
$$\mathcal{V}_{-4,2}(N):=\{n\in N\mid n\text{ is primitive},\ n^2=-4,\ \operatorname{div}_N(n)=2\}.$$
For $n\in\mathcal{V}_{-4,2}(N)$, define
$$s_n(x):=x+\frac{(x,n)}{2}n,\qquad\sigma_n:=[s_n]\in Q_N.$$
The divisibility assumption makes $s_n$ an integral isometry. It preserves $\mathcal{D}_N^\circ$, since it fixes pointwise a positive two-plane in $n^\perp$. Its fixed locus in $\mathcal{D}_N^\circ$ is the mirror
$$\mathcal{H}_n:=\{P\in\mathcal{D}_N^\circ\mid\underline{P}\subset n^\perp\},\qquad\mathcal{H}_n^\circ:=\mathcal{H}_n\cap\mathcal{U}_N.$$

\begin{lemma}[The retained mirrors]
\label{lem:retained-mirrors}
The space $\mathcal{H}_n^\circ$ is nonempty and path connected. At a general point of $\mathcal{H}_n^\circ$, the stabilizer in $Q_N$ is the order-two group generated by $\sigma_n$. In a normal holomorphic coordinate to the mirror, its nonidentity element acts by $z\mapsto-z$.
\end{lemma}

\begin{proof}
The orthogonal complement of $n$ has signature $(2,9)$, so $\mathcal{H}_n$ is a connected, contractible type-IV domain of complex dimension nine. No root mirror contains $\mathcal{H}_n$: such containment would force its root to be proportional to $n$, which is impossible because $n$ is primitive and has square $-4$.

The intersections with the root mirrors therefore form a locally finite family of proper complex hypersurfaces, with empty intersections omitted. Their complement is nonempty, and general position for paths shows that it is path connected. This complement is $\mathcal{H}_n^\circ$.

The pointwise stabilizer of $\mathcal{H}_n$ in the projective orthogonal group is $\{1,\sigma_n\}$. Indeed, an orthogonal transformation fixing all its period lines must act as a scalar on $n^\perp$ and independently as a scalar on $\mathbb{R}n$. The scalars are $\pm1$, and quotienting by the common sign leaves precisely the identity and the reflection.

Every other element of $Q_N$ consequently has a proper fixed locus on $\mathcal{H}_n$. Proper discontinuity makes these additional fixed loci locally finite, so they can be avoided. Finally, $s_n$ is the identity on $n^\perp$ and minus the identity on $\mathbb{R}n$, giving the stated action on the normal coordinate.
\end{proof}

\begin{proposition}[Kantor metric stabilizers]
\label{prop:kantor-metric-stabilizers}
For every $n\in\mathcal{V}_{-4,2}(N)$, there is an involution $\widetilde\kappa_n\in\operatorname{O}^{+,+}(\Lambda)^\tau$ whose restriction to $N$ is $s_n$. Its effective transformation $\gamma_n\in\Gamma_E^{+}$ fixes a point $x=(\ell,P)\in\mathcal{U}_E$, where $P\in\mathcal{H}_n^\circ$ has stabilizer $(Q_N)_P=\langle\sigma_n\rangle$. Thus $\gamma_n$ is realized by an involutive isometry of a smooth metric in the prescribed component $\mathcal{T}_E$, and
$$\rho_N(\gamma_n)=\sigma_n.$$
\end{proposition}

\begin{proof}
The class
$$a_n:=n/2+N\in A_N$$
is nonzero and satisfies $q_N(a_n)=-1\bmod2\mathbb{Z}$. For $x\in N^\vee$, the induced discriminant action is
$$\overline{s_n}(x+N)=(x+N)+\bigl((x,n)\bmod2\bigr)a_n.$$
This is the orthogonal transvection in $a_n$.

Set $a:=\varphi^{-1}(a_n)\in A_M$. Under $A_M\cong V=L/2L$, it is an anisotropic vector. The finite quadratic space $(V,q_2)$ is split, and its orthogonal group acts transitively on anisotropic vectors. This also follows directly from the model $q_2(x,y)=x\cdot y$ on $\mathbb{F}_2^5\oplus\mathbb{F}_2^5$: the subgroup acting by $(x,y)\mapsto(Ax,A^{-\mathsf{T}}y)$ is transitive on pairs with $x\cdot y=1$.

Start with a root in the $E_8(-1)$ summand of $L$. Surjectivity of reduction modulo two allows us to transport this root and the decomposition of $L$ to obtain
$$L=U_\nu\oplus E_\nu,\qquad U_\nu\cong U,\qquad E_\nu\cong E_8(-1),$$
with a root $\nu\in E_\nu$ satisfying
$$ (\nu,\nu)_L=-2,\qquad(\nu,\nu)_M=-4,\qquad\nu/2+M=a.$$
Here $L$ and $M=L(2)$ have the same underlying abelian group, and the two subscripts distinguish their bilinear forms. The required reduction homomorphism is the one used in Proposition~\ref{prop:sheet-preserving-arithmetic-action}.

Let $b_\nu$ be the Bertini involution associated with this decomposition, and put $F_\nu:=\nu^\perp\cap E_\nu\cong E_7(-1)$. Using reflection in $\nu$ on $L$, define
$$r_\nu(y):=y+(y,\nu)_L\nu,\qquad\kappa_\nu:=b_\nu r_\nu.$$
Then $\kappa_\nu$ is the identity on $U_\nu\oplus\mathbb{Z}\nu$ and minus the identity on $F_\nu$. This is the Kantor lattice involution described in \cite[Section~3]{Allcock2018}. Since $b_\nu$ is trivial modulo two, $\kappa_\nu$ induces the transvection in $a$ on $A_M$. The two discriminant actions therefore satisfy
$$\overline{s_n}\circ\varphi=\varphi\circ\overline{\kappa_\nu}.$$
Lemma~\ref{lem:gluing-extension} extends $(\kappa_\nu,s_n)$ to an involution $\widetilde\kappa_n$ of $\Lambda$ commuting with $\tau$. Its fixed subspaces in $M_{\mathbb{R}}$ and $N_{\mathbb{R}}$ have signatures $(1,2)$ and $(2,9)$, respectively. It follows that $\widetilde\kappa_n\in\operatorname{O}^{+,+}(\Lambda)^\tau$.

We compute the realization lattice for
$$K_n:=\langle\tau,\widetilde\kappa_n\rangle\cong(\mathbb{Z}/2\mathbb{Z})^2.$$
Its invariant positive three-planes are exactly $\ell\oplus\underline{P}$ with $\ell$ a positive line in $M_{\mathbb{R}}^{\kappa_\nu=+1}$ and $P\in\mathcal{H}_n$. The negative eigenspaces contain no positive directions, so the span of these planes has integral orthogonal complement
$$L_{K_n}=\Lambda\cap\bigl(F_\nu(2)_{\mathbb{R}}\oplus\mathbb{R}n\bigr).$$

Choose a root $r\in E_\nu$ with $(r,\nu)_L=-1$, and put $w:=r-\nu/2$. Such a root exists in the $E_8$ root system. We have
$$w\in F_\nu^\vee,\qquad(w,w)_L=-\frac32,\qquad w+\nu/2=r\in E_\nu.$$
Moreover, $2w\in F_\nu$, so $w\in M^\vee$, and $w+M=\nu/2+M=a$. The gluing condition therefore gives
$$z:=w+n/2\in\Lambda,\qquad z^2=-4.$$
This vector belongs to $L_{K_n}$.

Since $\operatorname{div}_N(n)=2$, we have
$$N^\vee\cap\mathbb{R}n=\tfrac12\mathbb{Z}n.$$
An integral multiple of $n$ has trivial class in $A_N$, forcing the $M$-component of a vector in $L_{K_n}$ to lie in $F_\nu(2)$. An odd half-integral multiple has class $a_n$, forcing that component to lie in $w+F_\nu(2)$. Hence the full realization lattice is the index-two extension
$$L_{K_n}=\bigl(F_\nu(2)\oplus\mathbb{Z}n\bigr)+\mathbb{Z}(w+n/2).$$

Both $n\in N$ and $\nu\in E_\nu(2)$ have square $-4$. Hence the map from $F_\nu(2)\oplus\mathbb{Z}n$ to $F_\nu(2)\oplus\mathbb{Z}\nu\subset E_\nu(2)$ that is the identity on $F_\nu$ and sends $n$ to $\nu$ is an isometry. It sends the extra generator to $w+\nu/2=r$. Pairing with $\nu$ modulo two gives
$$E_\nu=(F_\nu\oplus\mathbb{Z}\nu)+\mathbb{Z}r.$$
Consequently,
$$L_{K_n}\cong E_\nu(2)\cong E_8(-2).$$
In particular, $L_{K_n}$ contains no roots.

Proposition~\ref{prop:finite-group-realization} now gives an open dense set of $K_n$-invariant positive three-planes avoiding every K3 root wall. The fixed period domain projects openly and surjectively onto $\mathcal{H}_n$. We may therefore also require that $P$ lie in the open dense locus of $\mathcal{H}_n^\circ$ with stabilizer $\langle\sigma_n\rangle$ from Lemma~\ref{lem:retained-mirrors}.

Realization and descent give an involutive isometry in the prescribed metric component. Its period action is $\gamma_n$, and its projective restriction is $\sigma_n$, as required.
\end{proof}

We now record the resulting relations in the period action braid group. For $P\in\mathcal{H}_n^\circ$ and a path $\eta$ from $P_0$ to $P$, set
$$s_N(\sigma_n,P,\eta):=[\eta*(\sigma_n\eta)^{-1},\sigma_n]\in\mathcal{B}_N.$$
Let $\mathcal{R}_{-4}$ be the normal subgroup of $\mathcal{B}_N$ generated by all these classes. At a general point of the mirror, this class is the generator of its order-two isotropy group in the quotient orbifold.

Choose a Kantor lift $\gamma_n$ as above for each $n$. Let $\mathcal{N}_{\mathrm{Kan}}$ be the normal subgroup of $\mathcal{B}_E^{+}$ generated by its metric stabilizer classes $s(\gamma_n,x,\beta)$, allowing all such fixed points and paths. Define
$$\mathcal{N}_0:=\langle\mathcal{N}_{\mathrm{Bert}},\mathcal{N}_{\mathrm{Kan}}\rangle\subseteq\mathcal{N}_{\mathrm{stab}}^{+}.$$

\begin{corollary}[Kantor relations in the period braid group]
\label{cor:kantor-braid-relations}
The period projection satisfies
$$p_{\#}(\mathcal{N}_{\mathrm{Kan}})=\mathcal{R}_{-4},\qquad\mathcal{B}_E^{+}/\mathcal{N}_0\cong\mathcal{B}_N/\mathcal{R}_{-4}.$$
In particular, the ordinary fundamental group of $\Gamma_E^{+}\backslash\mathcal{T}_E$ is a quotient of $\mathcal{B}_N/\mathcal{R}_{-4}$.
\end{corollary}

\begin{proof}
For a Kantor metric stabilizer with $x=(\ell,P)$, equivariance gives
$$p_{\#}s(\gamma_n,x,\beta)=s_N(\sigma_n,P,p_N\circ\beta).$$
Thus every such image belongs to $\mathcal{R}_{-4}$.

For fixed $n$, these classes are conjugate: the mirror is path connected and fixed pointwise by $\sigma_n$, and changing a joining path conjugates its stabilizer class.

Proposition~\ref{prop:kantor-metric-stabilizers} supplies a metric stabilizer over at least one point of each mirror. Since $p_{\#}$ is surjective, the image of $\mathcal{N}_{\mathrm{Kan}}$ is normal in $\mathcal{B}_N$ and therefore contains all the generators of $\mathcal{R}_{-4}$. This proves the first equality.

The kernel of $p_{\#}$ is $\mathcal{N}_{\mathrm{Bert}}$, by Theorem~\ref{thm:bertini-kernel-removal}. The quotient isomorphism follows. Finally, $\mathcal{N}_0\subseteq\mathcal{N}_{\mathrm{stab}}^{+}$, so Armstrong's theorem gives the last assertion.
\end{proof}

\section{Modular forms and the residual period group}
\label{sec:residual-period-group}

\subsection{The Borcherds chart}
\label{subsec:borcherds-chart}

To compute $\mathcal{B}_N/\mathcal{R}_{-4}$, set
$$X_N:=Q_N\backslash\mathcal{D}_N^\circ,\qquad Y_N:=Q_N\backslash\mathcal{U}_N.$$
Consider the rescaled dual lattice
$$N':=N^\vee(2)\cong U\oplus U(2)\oplus E_8(-1)\cong 2U\oplus D_8(-1).$$
Here $D_8$ is positive definite; see \cite[Section~2]{GritsenkoHulek2016}. Duality and positive rescaling preserve the isometry group and period domain, giving
$$\operatorname{O}^{+}(N)=\operatorname{O}^{+}(N')\cong\operatorname{O}^{+}(2U\oplus D_8(-1)).$$

Let $\mathscr{A}_N$ be the graded ring of holomorphic modular forms with trivial character for $\operatorname{O}^{+}(N)$, graded by their usual automorphic weight:
$$\mathscr{A}_N:=\bigoplus_{k\geq 0}M_k(\operatorname{O}^{+}(N)).$$
We use the convention that a form of weight $k$, viewed on the affine cone over the period domain, satisfies $F(t\omega)=t^{-k}F(\omega)$. Since $-\operatorname{id}_N$ belongs to $\operatorname{O}^{+}(N)$, a nonzero form with trivial character has even weight. The theorem of Wang--Williams \cite[Theorem~5.2]{WangWilliams2020} states that $\mathscr{A}_N$ is a polynomial ring on eleven generators of weights
$$4,\ 6,\ 8,\ 8,\ 10,\ 10,\ 12,\ 12,\ 14,\ 16,\ 18.$$

Borcherds constructed a modular form $\Phi_4$ of weight four with a nontrivial quadratic character $\chi$ whose zero divisor on $\mathcal{D}_N^\circ$ is precisely the union of the $(-2)$-hyperplanes, each with multiplicity one \cite[Theorem~3.2 and Lemmas~3.6,~3.8]{Borcherds1996}. In particular,
$$\operatorname{div}_{\mathcal{D}_N^\circ}(\Phi_4)=\sum_{\delta\in\operatorname{Roots}(N)/\{\pm1\}}D_\delta.$$
Its square $q:=\Phi_4^2\in M_8(\operatorname{O}^+(N))$ has trivial character.

\begin{lemma}
\label{lem:borcherds-generator}
The polynomial generators of $\mathscr{A}_N$ can be chosen so that one of the weight-eight generators is $q$. More precisely, there are algebraically independent generators such that
$$\mathscr{A}_N=\mathbb{C}[f_4,f_6,q,f_8,f_{10,1},f_{10,2},f_{12,1},f_{12,2},f_{14},f_{16},f_{18}],$$
where each subscript records the automorphic weight.
\end{lemma}

\begin{proof}
Start with any homogeneous polynomial generators of the weights listed above, and denote the two weight-eight generators by $h_8$ and $h_8'$. The weight-eight part of the ring is spanned by $f_4^2,h_8,h_8'$, so
$$q=a f_4^2+b h_8+c h_8'$$
for constants $a,b,c\in\mathbb{C}$. If $b=c=0$, then $a\neq0$ and the identity $\Phi_4^2=a f_4^2$ implies that $\Phi_4$ is a constant multiple of $f_4$ on the connected affine cone over $\mathcal{D}_N^\circ$. This contradicts the nontriviality of the character of $\Phi_4$, since $f_4$ has trivial character. Hence $(b,c)\neq(0,0)$.

Choose $f_8$ independent of $b h_8+c h_8'$ in their two-dimensional span. Replacing $h_8,h_8'$ by $q,f_8$ is an invertible graded polynomial change of generators.
\end{proof}

Let $X_N^{\mathrm{BB}}$ denote the Baily--Borel compactification of $X_N$. Its description by modular forms gives
$$X_N^{\mathrm{BB}}\cong\operatorname{Proj}(\mathscr{A}_N).$$
We now divide all degrees in this ring by two. This regrading leaves its coarse projective spectrum unchanged and yields
$$X_N^{\mathrm{BB}}\cong\mathbb{P}(2,3,4,4,5,5,6,6,7,8,9).$$
In the normalized grading, $q$ has degree four. The open chart $D_+(q)$ is the complement of its zero locus in $X_N^{\mathrm{BB}}$. Writing $\partial X_N^{\mathrm{BB}}:=X_N^{\mathrm{BB}}\setminus X_N$, the divisor description gives
$$Y_N=D_+(q)\setminus\partial X_N^{\mathrm{BB}}.$$
Every Baily--Borel boundary stratum for this type IV domain has complex dimension zero or one; in particular, $\dim_{\mathbb{C}}\partial X_N^{\mathrm{BB}}\leq1$. The boundary strata in the Enriques case are described explicitly in \cite[Section~1]{Borcherds1996}.

Let $x_1,\ldots,x_{10}$ correspond, in order, to the generators other than $q$ in Lemma~\ref{lem:borcherds-generator}. Their normalized degrees are
$$d=(d_1,\ldots,d_{10})=(2,3,4,5,5,6,6,7,8,9).$$
The group $\mu_4:=\{\zeta\in\mathbb{C}^*\mid\zeta^4=1\}$ acts on $\mathbb{A}^{10}_{\mathbb{C}}$ by
$$\zeta\cdot(x_1,\ldots,x_{10})=(\zeta^{d_1}x_1,\ldots,\zeta^{d_{10}}x_{10}).$$

\begin{proposition}
\label{prop:borcherds-affine-chart}
There is a finite quotient map
$$\pi_q:\mathbb{A}^{10}_{\mathbb{C}}\longrightarrow D_+(q)\cong\mathbb{A}^{10}_{\mathbb{C}}/\mu_4.$$
Define
$$B_\partial:=\pi_q^{-1}\bigl(\partial X_N^{\mathrm{BB}}\cap D_+(q)\bigr),\qquad Z_N:=\mathbb{A}^{10}_{\mathbb{C}}\setminus B_\partial.$$
Then $B_\partial$ is a $\mu_4$-invariant closed algebraic subset of complex dimension at most one, and
$$Y_N\cong Z_N/\mu_4.$$
The space $Z_N$ is connected and simply connected in the analytic topology. The $\mu_4$-action is effective, and the fixed locus of every nonidentity element has complex codimension at least five in $Z_N$.
\end{proposition}

\begin{proof}
In the normalized grading, the weighted scalar action on the affine coordinates of $\operatorname{Spec}(\mathscr{A}_N)$ is
$$t\cdot(q,x_1,\ldots,x_{10})=(t^4q,t^{d_1}x_1,\ldots,t^{d_{10}}x_{10}).$$
On the locus $q\neq0$, every scalar orbit meets the slice $q=1$, and two points of that slice lie in the same orbit precisely when they differ by the indicated $\mu_4$-action. Equivalently, restriction to $q=1$ identifies the degree-zero localization with the invariant ring:
$$\bigl(\mathscr{A}_N[q^{-1}]\bigr)_0\cong\mathbb{C}[x_1,\ldots,x_{10}]^{\mu_4}.$$
Taking spectra proves the description of $D_+(q)$ and defines the finite map $\pi_q$.

The boundary is closed in $X_N^{\mathrm{BB}}$, so $B_\partial$ is closed and algebraic. Since $\pi_q$ is finite and the boundary has dimension at most one, the same dimension bound holds for $B_\partial$. Removing its inverse image from the quotient presentation gives $Y_N\cong Z_N/\mu_4$.

Choose a Whitney stratification of $B_\partial$. Every stratum has real codimension at least eighteen in $\mathbb{C}^{10}$. Relative transversality therefore allows any path with endpoints in $Z_N$ to be moved away from $B_\partial$. It also allows any disk in $\mathbb{C}^{10}$ whose boundary lies in $Z_N$ to be moved away from $B_\partial$ while keeping its boundary fixed. Since $\mathbb{C}^{10}$ is contractible, these two observations prove that $Z_N$ is connected and simply connected.

Finally, reduce the weights modulo four:
$$d\equiv(2,3,0,1,1,2,2,3,0,1)\pmod4.$$
A primitive fourth root fixes exactly the coordinate subspace spanned by $x_3,x_9$, of dimension two. Its square fixes the coordinate subspace spanned by $x_1,x_3,x_6,x_7,x_9$, of dimension five. The same dimension-two description applies to the other primitive fourth root. These subspaces meet $Z_N$, since none can be contained in $B_\partial$. Thus the nonidentity fixed loci have codimensions eight or five. The action is effective because, for example, the weight of $x_2$ is relatively prime to four.
\end{proof}

\begin{remark}
\label{rem:borcherds-chart-stacks}
Dividing the weights by two removes the generic ineffective $\mu_2$ of the original weighted quotient stack. The resulting stack over $Y_N$ is $[Z_N/\mu_4]$. Its coarse presentation alone does not identify $\pi_1(Y_N)$ with $\mu_4$, because the action has fixed points. We next compare this stack with the period quotient stack.
\end{remark}

\subsection{The canonical stack and the residual group}
\label{subsec:canonical-stack}

Let
$$\mathscr{P}_N:=[\mathcal{U}_N/Q_N].$$
Write $\pi_1^{\mathrm{orb}}$ for the fundamental group of the analytic stack, or its classifying space. For this discrete global quotient it is the action braid group:
$$\pi_1^{\mathrm{orb}}(\mathscr{P}_N,P_0)\cong\mathcal{B}_N.$$
Basepoints in the subsequent comparisons are chosen compatibly; changing connecting paths only conjugates the homomorphisms and stabilizer images involved.

The orbifold $\mathscr{P}_N$ has an algebraic model. Choose a normal neat subgroup $Q_N'\subset Q_N$ of finite index. It acts freely, and $S_N:=Q_N'\backslash\mathcal{U}_N$ is a smooth quasi-projective variety: it is the complement of the pullback of the Borcherds divisor in the corresponding arithmetic quotient. The finite group $F_N:=Q_N/Q_N'$ acts algebraically on $S_N$, and
$$\mathscr{P}_N\cong[S_N/F_N]^{\mathrm{an}}.$$
Thus $\mathscr{P}_N$ is the analytification of a smooth separated Deligne--Mumford stack of finite type over $\mathbb{C}$. Its generic stabilizer is trivial, because the projective action of $Q_N$ is effective. Its coarse space is $Y_N$.

\begin{lemma}[Divisorial stabilizers]
\label{lem:divisorial-period-stabilizers}
The divisorial fixed loci of nonidentity elements of $Q_N$ in $\mathcal{U}_N$ are precisely the retained mirrors $\mathcal{H}_n^\circ$, with $n\in\mathcal{V}_{-4,2}(N)$. Their generic stabilizers have order two and are generated by $\sigma_n$.

Consequently, under $\pi_1^{\mathrm{orb}}(\mathscr{P}_N)\cong\mathcal{B}_N$, the normal subgroup generated by the images of the generic divisorial stabilizers is exactly $\mathcal{R}_{-4}$.
\end{lemma}

\begin{proof}
We recall the reflection calculation underlying the description of the ramification divisor in \cite[Section~5]{GritsenkoHulek2016}. Suppose that a nonidentity projective isometry $[g]\in Q_N$ fixes a divisor in $\mathcal{D}_N^\circ$. Since point stabilizers are finite, a representative $g\in\operatorname{O}^{+}(N)$ has finite order. Its projective fixed locus is the union of the intersections of the period domain with its projectivized complex eigenspaces.

An eigenspace with eigenvalue $\lambda\notin\{1,-1\}$ is totally isotropic for the complex bilinear form. Its dimension is therefore at most six, so its projectivization cannot contain a divisor of the ten-dimensional period domain. An eigenspace with eigenvalue $1$ or $-1$ is nondegenerate, and its projectivized quadric has dimension two less than the eigenspace. To contain the fixed divisor, it must therefore have dimension eleven. After changing the sign of $g$, we conclude that $g$ is the identity on a rational hyperplane and minus the identity on its orthogonal line. That line is negative, since the hyperplane contains a positive two-plane. Hence $g=s_v$ for a primitive negative vector $v\in N$.

Put $d:=\operatorname{div}_N(v)$. The class of $v/d$ in $A_N$ has order $d$. Since $A_N$ has exponent two, $d\in\{1,2\}$. Integrality of
$$s_v(x)=x-\frac{2(x,v)}{v^2}v$$
requires $v^2\mid2d$. As $N$ is even, the possibilities are
$$\bigl(v^2,d\bigr)=(-2,1),\quad(-2,2),\quad(-4,2).$$
The discriminant quadratic form of $N=U\oplus U(2)\oplus E_8(-2)$ takes values in $\mathbb{Z}/2\mathbb{Z}$. The middle possibility would instead give
$$q_N(v/2+N)=-\tfrac12\bmod2\mathbb{Z},$$
so it is excluded. The two reflective types are therefore $(-2,1)$ and $(-4,2)$.

The $(-2)$-mirrors are absent from $\mathcal{U}_N$. The remaining mirrors are exactly $\mathcal{H}_n^\circ$ with $n\in\mathcal{V}_{-4,2}(N)$, and Lemma~\ref{lem:retained-mirrors} identifies their generic stabilizers as $\langle\sigma_n\rangle$.

The image of this inertia generator in the action fundamental group is the class $s_N(\sigma_n,P,\eta)$ defined in Section~\ref{subsec:kantor-stabilizers}. Allowing arbitrary points of a retained mirror gives the same normal subgroup as using general points: the mirror is path connected, and the reflection fixes it pointwise, so the resulting classes differ by conjugation. This is precisely the defining normal subgroup $\mathcal{R}_{-4}$.
\end{proof}

For a complex variety $Y$ with quotient singularities, its canonical stack $Y^{\mathrm{can}}$ is the smooth Deligne--Mumford stack with coarse space $Y$ whose coarse morphism is an isomorphism away from a subset of complex codimension at least two. The characterization and universal property recalled in \cite[Section~2]{GeraschenkoSatriano2017} determine it up to isomorphism over $Y$.

\begin{lemma}[Removing divisorial inertia]
\label{lem:canonical-stack-fundamental-group}
Let $\mathscr{X}$ be a connected smooth separated Deligne--Mumford stack of finite type over $\mathbb{C}$ with trivial generic stabilizer, and let $Y$ be its coarse space. The canonical comparison morphism
$$c:\mathscr{X}\longrightarrow Y^{\mathrm{can}}$$
induces a surjective homomorphism on topological orbifold fundamental groups. Its kernel is the normal subgroup $\mathcal{I}_{\mathrm{div}}$ generated by the images of the inertia groups at general points of the divisorial stabilizer loci. Thus
$$\pi_1^{\mathrm{orb}}(\mathscr{X})/\mathcal{I}_{\mathrm{div}}\cong\pi_1^{\mathrm{orb}}(Y^{\mathrm{can}}).$$
\end{lemma}

\begin{proof}
Let $D_1,\ldots,D_r$ be the irreducible components of the branch divisor in $Y$, and let $e_i$ be the generic inertia order along $D_i$. The structure theorem of Geraschenko--Satriano \cite[Theorem~1]{GeraschenkoSatriano2017} supplies $c$ and expresses $\mathscr{X}$ as the canonical stack of the iterated root stack over $Y^{\mathrm{can}}$ with orders $e_i$ along these divisors.

Remove from $Y$ its singular locus, the singular loci and pairwise intersections of the $D_i$, and the remaining loci where stabilizers differ from their generic divisorial values. These sets are contained in a closed algebraic subset $C\subset Y$ of complex codimension at least two. Write
$$V:=Y\setminus C,\qquad W:=V\setminus\bigcup_{i=1}^rD_i.$$
Over $V$, the canonical stack is just $V$, and $\mathscr{X}_V:=\mathscr{X}\times_YV$ is the iterated root stack along the disjoint smooth divisors $D_i\cap V$.

Deleting the inverse image of $C$ from either smooth stack does not change its orbifold fundamental group. Indeed, in local finite quotient charts the deleted strata have real codimension at least four; relative general position for paths and homotopies, together with orbifold van Kampen, gives the assertion. Therefore
$$\pi_1^{\mathrm{orb}}(\mathscr{X}_V)\cong\pi_1^{\mathrm{orb}}(\mathscr{X}),\qquad\pi_1(V)\cong\pi_1^{\mathrm{orb}}(Y^{\mathrm{can}}).$$

Choose based meridians $m_i\in\pi_1(W)$ around the divisors $D_i\cap V$. Locally, a root construction has the normal-coordinate model
$$[\Delta_z/\mu_{e_i}]\longrightarrow\Delta_t,\qquad t=z^{e_i},$$
with the group acting trivially in the directions tangent to the divisor. A meridian in the punctured $t$-disk becomes the local inertia generator; its $e_i$th power bounds a disk in the uniformizing chart. Orbifold van Kampen consequently gives
$$\pi_1^{\mathrm{orb}}(\mathscr{X}_V)\cong\pi_1(W)/\langle\!\langle m_1^{e_1},\ldots,m_r^{e_r}\rangle\!\rangle.$$
Ordinary van Kampen gives
$$\pi_1(V)\cong\pi_1(W)/\langle\!\langle m_1,\ldots,m_r\rangle\!\rangle.$$
The comparison map is induced by the identity on $W$. It is therefore surjective, and its kernel is normally generated by the images of the $m_i$. Each such image is the image of a generator of the corresponding generic inertia group. This proves the assertion.
\end{proof}

\begin{theorem}[The residual period group]
\label{thm:residual-period-group}
There are isomorphisms
$$\mathcal{B}_N/\mathcal{R}_{-4}\cong\pi_1^{\mathrm{orb}}(Y_N^{\mathrm{can}})\cong\mu_4.$$
Consequently,
$$\mathcal{B}_E^{+}/\mathcal{N}_0\cong\mu_4,$$
and the ordinary fundamental group of $\Gamma_E^{+}\backslash\mathcal{T}_E$ is a quotient of $\mu_4$.
\end{theorem}

\begin{proof}
The algebraic model of $\mathscr{P}_N$ satisfies the hypotheses of Lemma~\ref{lem:canonical-stack-fundamental-group}. Lemma~\ref{lem:divisorial-period-stabilizers} identifies its divisorial inertia subgroup with $\mathcal{R}_{-4}$. Applying the comparison lemma gives
$$\mathcal{B}_N/\mathcal{R}_{-4}\cong\pi_1^{\mathrm{orb}}(Y_N^{\mathrm{can}}).$$

By Proposition~\ref{prop:borcherds-affine-chart}, the smooth variety $Z_N$ is simply connected and $Y_N\cong Z_N/\mu_4$. The action is effective, and every nonidentity fixed locus has complex codimension at least five. Thus the coarse morphism $[Z_N/\mu_4]\to Y_N$ is an isomorphism away from a subset of codimension at least two. The characterization of the canonical stack therefore gives
$$Y_N^{\mathrm{can}}\cong[Z_N/\mu_4].$$
Equivalently, this follows from the absence of pseudoreflections in the local quotient charts; see \cite[Section~4]{GeraschenkoSatriano2017}.

For $z_0\in Z_N$, the action fundamental group has the exact sequence
$$1\longrightarrow\pi_1(Z_N,z_0)\longrightarrow\pi_1^{\mathrm{orb}}([Z_N/\mu_4],z_0)\longrightarrow\mu_4\longrightarrow1.$$
The left-hand group is trivial. Hence $\pi_1^{\mathrm{orb}}([Z_N/\mu_4])\cong\mu_4$, proving the first assertion. The remaining statements follow from Corollary~\ref{cor:kantor-braid-relations}.
\end{proof}

Let $\vartheta_N:\mathcal{B}_N\to\mu_4$ be this quotient homomorphism, with target identified using the normalized-coordinate action. The theorem states that
$$\ker\vartheta_N=\mathcal{R}_{-4},\qquad\ker(\vartheta_N\circ p_{\#})=\mathcal{N}_0.$$
It remains to find a smooth metric stabilizer whose image generates $\mu_4$.

\subsection{The Ohashi stabilizer and the sheet-preserving quotient}
\label{subsec:ohashi-stabilizer}

For $h\in Q_N$ fixing $P\in\mathcal{U}_N$ and a path $\eta$ from $P_0$ to $P$, extend the notation of Section~\ref{subsec:kantor-stabilizers} by setting
$$s_N(h,P,\eta):=[\eta*(h\eta)^{-1},h]\in\mathcal{B}_N.$$

\begin{lemma}[The residual character on stabilizers]
\label{lem:residual-stabilizer-character}
Let $u\in\operatorname{O}^{+}(N)$ fix $P\in\mathcal{U}_N$ projectively. Choose a nonzero period vector $\omega$ and write $u\omega=a\omega$.
Then $a^8=1$, and, with the normalized-coordinate action of Section~\ref{subsec:borcherds-chart},
$$\vartheta_N\bigl(s_N([u],P,\eta)\bigr)=a^2\in\mu_4.$$
The right-hand side is unchanged when $u$ is replaced by $-u$, and the image is independent of the choice of $\eta$.
\end{lemma}

\begin{proof}
For the affine period cone and the normalized slice, put
$$\mathcal{C}_N:=\{\omega\in N_{\mathbb{C}}\setminus\{0\}\mid[\omega]\in\mathcal{U}_N\}.$$
$$\mathcal{C}_{N,q}:=\{\omega\in\mathcal{C}_N\mid q(\omega)=1\}.$$
Since $q(t\omega)=t^{-8}q(\omega)$ and $q$ is nowhere zero on $\mathcal{C}_N$, the projection $\mathcal{C}_{N,q}\to\mathcal{U}_N$ is a principal $\mu_8$-cover, with $\mu_8$ acting by scalar multiplication.

The group acting on this cover that accounts for the effective projective action is
$$\widehat{Q}_{N,q}:=(\mu_8\times\operatorname{O}^{+}(N))/\langle(-1,-\operatorname{id}_N)\rangle.$$
The class of $(t,u)$ acts by $\omega\mapsto tu\omega$. Scalar multiplication gives an exact sequence
$$1\longrightarrow\mu_8\longrightarrow\widehat{Q}_{N,q}\longrightarrow Q_N\longrightarrow1,$$
and quotienting first by the freely acting scalar subgroup gives
$$[\mathcal{C}_{N,q}/\widehat{Q}_{N,q}]\cong[\mathcal{U}_N/Q_N]=\mathscr{P}_N.$$

Denote the ten modular generators other than $q$ by $F_1,\ldots,F_{10}$, ordered as in Section~\ref{subsec:borcherds-chart}. Their weights are $2d_1,\ldots,2d_{10}$. Evaluation defines
$$\operatorname{ev}_q:\mathcal{C}_{N,q}\longrightarrow Z_N,\qquad\omega\longmapsto(F_1(\omega),\ldots,F_{10}(\omega)).$$
The image lies in $Z_N$ because the period belongs to the Baily--Borel interior and $q=1$. Since $F_j(tu\omega)=t^{-2d_j}F_j(\omega)$, evaluation is equivariant for
$$\rho_q:\widehat{Q}_{N,q}\longrightarrow\mu_4,\qquad[(t,u)]\longmapsto t^{-2}.$$
This is well defined on the displayed quotient. It induces a morphism
$$\mathscr{P}_N\longrightarrow[Z_N/\mu_4]\cong Y_N^{\mathrm{can}}$$
over the identity of the coarse space $Y_N$.

This morphism agrees with the canonical comparison over the inverse image of the smooth locus of $Y_N$, where the target has trivial stabilizers. The agreement extends uniquely: in a smooth local chart of the source, the isomorphism space between the two morphisms is finite, and the closure of the section over the dense open subset is finite and birational over a normal space. It is therefore an isomorphism. Consequently, this evaluation construction induces the homomorphism $\vartheta_N$ used in Section~\ref{subsec:canonical-stack}.

Rescale the chosen $\omega$ so that $q(\omega)=1$. Invariance of $q$ gives
$$1=q(u\omega)=q(a\omega)=a^{-8},$$
so $a\in\mu_8$. The element $[u]\in(Q_N)_P$ is represented in the stabilizer of $\omega\in\mathcal{C}_{N,q}$ by $[(a^{-1},u)]$. Its image under $\rho_q$ is
$$\rho_q([(a^{-1},u)])=a^2.$$
Passing from this local stabilizer to the based action fundamental group gives the asserted formula. Replacing $u$ by $-u$ changes $a$ to $-a$, leaving its square unchanged. Changing the joining path conjugates the stabilizer class, and conjugation has no effect in the abelian group $\mu_4$.
\end{proof}

\begin{proposition}[Ohashi's generator]
\label{prop:ohashi-generator}
In every prescribed metric Teichm\"uller component $\mathcal{T}_E$, there is a smooth unit-volume Ricci-flat metric with a sheet-preserving finite-order isometry whose stabilizer class maps to a generator of $\mu_4$ under $\vartheta_N\circ p_{\#}$.
\end{proposition}

\begin{proof}
Ohashi constructs a complex Enriques surface $S$ with an automorphism $\sigma$ of order eight satisfying
$$\sigma^*\Omega_S=-i\Omega_S,\qquad0\ne\Omega_S\in H^0(S,K_S^{\otimes2});$$
see \cite[Example~1.2]{Ohashi2015}. Identify the underlying smooth manifold of $S$ with $E$ and its canonical cover with the fixed cover $p:X\to E$. We use the compatible marking $\alpha$ fixed in Section~\ref{subsec:enriques-lattices}. Choose a holomorphic lift $\widetilde\sigma$ and a nonzero holomorphic two-form $\omega_X$ on $X$, writing
$$\widetilde\sigma^*\omega_X=\lambda\omega_X.$$
Since $p^*\Omega_S$ is a nonzero multiple of $\omega_X^{\otimes2}$, equivariance gives $\lambda^2=-i$. Thus $\lambda$ is a primitive eighth root of unity. The other lift changes $\lambda$ to $-\lambda$.

To obtain an invariant metric, choose a K\"ahler class $\kappa_0$ on $S$ and average it:
$$\kappa:=\frac18\sum_{j=0}^{7}(\sigma^j)^*\kappa_0.$$
This class is K\"ahler and $\sigma$-invariant. Lemma~\ref{lem:averaging-kahler-class}, using the uniqueness in Yau's theorem, produces a smooth unit-volume Ricci-flat K\"ahler metric $g_0$ invariant under $\sigma$.

With the fixed marking and the inverse-pullback convention, set
$$A_\sigma:=\alpha\circ(\widetilde\sigma^{-1})^*\circ\alpha^{-1},\qquad u:=A_\sigma|_N.$$
Its action on $\alpha(\omega_X)$ has eigenvalue $\lambda^{-1}$. Since $\sigma$ preserves the averaged K\"ahler class and is holomorphic, $A_\sigma$ preserves the invariant positive line and the orientation of the anti-invariant positive two-plane. Its period transformation is therefore sheet preserving.

Apply Corollary~\ref{cor:transfer-finite-actions} to the finite isometry group $\langle\sigma\rangle$ of $(E,g_0)$. It gives an isomorphic finite group acting isometrically on a metric representing a point of the prescribed component $\mathcal{T}_E$, at the same marked metric period. The full marked action on the K3 cover is retained. In particular, the element corresponding to $A_\sigma$ gives $\gamma_\sigma\in\Gamma_E^{+}$ fixing a smooth metric period $x=(\ell,P)$, with
$$\rho_N(\gamma_\sigma)=[u].$$

The period line representing $P\in\mathcal{D}_N^\circ$ is either $\mathbb{C}\alpha(\omega_X)$ or its complex conjugate, according to the chosen component of the oriented period domain. The eigenvalue $a$ of $u$ on this line consequently satisfies
$$a^2\in\{\lambda^{-2},\overline{\lambda}^{\,-2}\}=\{i,-i\}.$$
Choose a path $\beta$ in $\mathcal{U}_E$ from $x_0$ to $x$. The period projection sends its metric stabilizer class to
$$p_{\#}\bigl(s(\gamma_\sigma,x,\beta)\bigr)=s_N([u],P,p_N\circ\beta).$$
By Lemma~\ref{lem:residual-stabilizer-character},
$$\vartheta_N\bigl(p_{\#}(s(\gamma_\sigma,x,\beta))\bigr)=a^2\in\{i,-i\}.$$
Both values generate $\mu_4$. This proves the proposition in the prescribed component.
\end{proof}

\begin{corollary}[The sheet-preserving coarse quotient]
\label{cor:sheet-preserving-simply-connected}
The coarse quotient by the sheet-preserving actual action is simply connected:
$$\pi_1(\Gamma_E^{+}\backslash\mathcal{T}_E)=\{1\}.$$
Equivalently,
$$\mathcal{B}_E^{+}=\mathcal{N}_{\mathrm{stab}}^{+}.$$
\end{corollary}

\begin{proof}
Theorem~\ref{thm:residual-period-group} gives
$$\mathcal{B}_E^{+}/\mathcal{N}_0\cong\mu_4,\qquad\mathcal{N}_0\subseteq\mathcal{N}_{\mathrm{stab}}^{+}.$$
Proposition~\ref{prop:ohashi-generator} supplies an element of $\mathcal{N}_{\mathrm{stab}}^{+}$ whose image generates this quotient. Hence
$$\mathcal{B}_E^{+}/\mathcal{N}_{\mathrm{stab}}^{+}=\{1\}.$$
The action-group form of Armstrong's theorem from Section~\ref{subsec:action-braid-groups} identifies this group with the ordinary fundamental group of $\Gamma_E^{+}\backslash\mathcal{T}_E$. The quotient is path connected because $\mathcal{T}_E$ is path connected, completing the proof.
\end{proof}

\section{The full quotient and proof of the main theorem}
\label{sec:proof-main}

\subsection{A sheet-reversing metric stabilizer}
\label{subsec:sheet-reversing-stabilizer}

\begin{proposition}[A realizable sheet reverser]
\label{prop:sheet-reversing-stabilizer}
Some smooth metric in $\mathcal{T}_E$ has an involutive stabilizer $\gamma_{\mathrm{rev}}\in\Gamma_E^{\mathrm{act}}$ with $\varepsilon_E(\gamma_{\mathrm{rev}})=1$. Consequently,
$$\Gamma_E^{\mathrm{act}}/\Gamma_E^+\cong\mathbb{Z}/2\mathbb{Z}.$$
\end{proposition}

\begin{proof}
In the unimodular hyperbolic summand of $N$, choose $e,f$ with $e^2=f^2=0$ and $(e,f)=1$. Set $v=e+f$. Since $v^2=2$, its reflection is the integral involution $s_v(x)=x-(x,v)v$. It acts trivially on $A_N$, and $-\operatorname{id}_M$ acts trivially on the $2$-elementary group $A_M$. Thus $(-\operatorname{id}_M,s_v)$ extends to an involution $a$ of $\Lambda$ commuting with $\tau$.

The group $K:=\langle\tau,a\rangle\cong(\mathbb{Z}/2\mathbb{Z})^2$ acts diagonally on
$$\Lambda_{\mathbb{R}}=M_{\mathbb{R}}\oplus\mathbb{R}v\oplus(v^\perp\cap N_{\mathbb{R}}),\qquad\tau=(+1,-1,-1),\quad a=(-1,-1,+1).$$
The summands have signatures $(1,9)$, $(1,0)$, and $(1,10)$. For positive lines $\ell$ in the first and $m$ in the third, the plane $\ell\oplus\mathbb{R}v\oplus m$ is $K$-invariant, with both generators of determinant $+1$. Hence $K\subset\operatorname{O}^+(\Lambda)$.

These planes span all of $\Lambda_{\mathbb{R}}$, so $W_K=\Lambda_{\mathbb{R}}$ and $L_K=0$. Proposition~\ref{prop:finite-group-realization} therefore realizes $K/\langle\tau\rangle$ on a smooth metric in the prescribed component. The class of $a$ gives $\gamma_{\mathrm{rev}}$. Reflection in the positive vector $v$ reverses the orientation of the positive plane $\mathbb{R}v\oplus m$, so $\varepsilon_E(\gamma_{\mathrm{rev}})=1$. At its fixed metric period, orientation correction gives $c\,s_vP=P$, as required by our action convention.
\end{proof}

\subsection{Completion of the proof}
\label{subsec:completion-main-proof}

\begin{proof}[Proof of Theorem~\ref{thm:main}]
Choose $\mathcal{T}_E$ mapping to the given coarse component $\mathcal{M}_0$. The sheet character and endpoint map give
$$1\longrightarrow\mathcal{B}_E^+\longrightarrow\mathcal{B}_E\xrightarrow{\varepsilon_E\circ e_E}\mathbb{Z}/2\mathbb{Z}\longrightarrow1.$$
By Corollary~\ref{cor:sheet-preserving-simply-connected},
$$\mathcal{B}_E^+=\mathcal{N}_{\mathrm{stab}}^+\subseteq\mathcal{N}_{\mathrm{stab}}.$$
A stabilizer braid for $\gamma_{\mathrm{rev}}$ lies in $\mathcal{N}_{\mathrm{stab}}$ and maps to the nontrivial element of the last quotient. Hence $\mathcal{N}_{\mathrm{stab}}=\mathcal{B}_E$, and Corollary~\ref{cor:coarse-braid-quotient} gives $\pi_1(\mathcal{M}_0)=1$. The component is path connected as a quotient of $\mathcal{T}_E$, and abelianization gives $H_1(\mathcal{M}_0;\mathbb{Z})=0$.

Finally, let $\mathcal{R}^{\mathrm{Ric}=0}(E)$ and $\mathcal{M}^{\mathrm{Ric}=0}(E)$ denote the corresponding spaces without volume normalization. In dimension four, $\operatorname{Vol}_{rg}(E)=r^2\operatorname{Vol}_g(E)$. Thus, with $v(g):=\operatorname{Vol}_g(E)$,
$$g\longmapsto\bigl(v(g)^{-1/2}g,v(g)\bigr)$$
is an equivariant homeomorphism $\mathcal{R}^{\mathrm{Ric}=0}(E)\cong\mathcal{R}^{\mathrm{Ric}=0}_{\mathrm{unit}}(E)\times\mathbb{R}_{>0}$, with inverse $(h,v)\mapsto v^{1/2}h$ and trivial action on the second factor. Since orbit maps are open, it descends to
$$\mathcal{M}^{\mathrm{Ric}=0}(E)\cong\mathcal{M}^{\mathrm{Ric}=0}_{\mathrm{unit}}(E)\times\mathbb{R}_{>0}.$$
The second factor is contractible, completing the proof.
\end{proof}

\section*{Acknowledgements}

The author is grateful to his advisor, Professor Dr.~Wilderich Tuschmann, for his supervision, guidance, and support throughout this work.

The author acknowledges financial support from the Research Training Group RTG~2229 \textquotedblleft Asymptotic Invariants and Limits of Groups and Spaces\textquotedblright, funded by the Deutsche Forschungsgemeinschaft (DFG, German Research Foundation).

The author used GPT~5.6~Sol to assist with the mathematical reasoning in Sections~\ref{sec:metric-stabilizers} and~\ref{sec:residual-period-group}. GPT~6~Astra was used to generate draft text during the preparation of the manuscript.
\bibliographystyle{amsplain}
\bibliography{references}

\end{document}